\documentclass{article}
\usepackage{amsthm}
\usepackage[sumlimits]{amsmath}
\usepackage{amsfonts}
\usepackage{color}
\usepackage{graphicx}
\DeclareMathOperator{\ad}{ad}
\DeclareMathOperator{\Vol}{Vol}
\DeclareMathOperator{\diff}{d}
\DeclareMathOperator{\Emb}{Emb}

\DeclareMathOperator{\Id}{Id}
\DeclareMathOperator{\Tr}{Tr}

\newtheorem{theorem}{Theorem}
\newtheorem{lemma}[theorem]{Lemma}
\newtheorem{corollary}[theorem]{Corollary}
\newtheorem{definition}[theorem]{Definition}
\newtheorem{example}[theorem]{Example}

\newtheorem{remark}[theorem]{Remark}
\def\MM#1{\boldsymbol{#1}}
\newcommand{\pp}[2]{\frac{\partial #1}{\partial #2}} 
\newcommand{\dede}[2]{\frac{\delta #1}{\delta #2}}
\newcommand{\dd}[2]{\frac{\diff#1}{\diff #2}} 

\newcommand{\eqnref}[1]{(\ref{#1})}

\newcommand{\logapp}{\overline{\log}}
\newcommand{\expapp}{\overline{\exp}}
\newcommand{\bfi}[1]{{\bfseries\itshape #1}}

\begin{document}

\title{Continuous and discrete Clebsch variational principles}
\author{C. J. Cotter and D. D. Holm}
\maketitle
\begin{abstract}
  The Clebsch method provides a unifying approach for deriving
  variational principles for continuous and discrete dynamical systems
  where elements of a vector space are used to control dynamics on the
  cotangent bundle of a Lie group \emph{via} a velocity map. This
  paper proves a reduction theorem which states that the canonical
  variables on the Lie group can be eliminated, if and only if the
  velocity map is a Lie algebra action, thereby producing the
  Euler-Poincar\'e (EP) equation for the vector space variables. In
  this case, the map from the canonical variables on the Lie group to
  the vector space is the standard momentum map defined using the
  diamond operator. We apply the Clebsch method in examples of the
  rotating rigid body and the incompressible Euler equations. Along
  the way, we explain how singular solutions of the EP equation for
  the diffeomorphism group (EPDiff) arise as momentum maps in the
  Clebsch approach. In the case of finite dimensional Lie groups, the
  Clebsch variational principle is discretised to produce a
  variational integrator for the dynamical system. We obtain a
  discrete map from which the variables on the cotangent bundle of a
  Lie group may be eliminated to produce a discrete EP equation for
  elements of the vector space. We give an integrator for the rotating
  rigid body as an example. We also briefly discuss how to discretise
  infinite-dimensional Clebsch systems, so as to produce conservative
  numerical methods for fluid dynamics.
\end{abstract}



\section{Introduction}

We are dealing with variational principles defined by an action (or cost, for the optimal control problem) 
\begin{equation}
S = \int l [\MM{\xi}(t)]\,dt 
\,,
\label{lag-cost}
\end{equation}
whose Lagrangian (or cost functional) $l:\,V\mapsto \mathbb{R}$ is defined on vectors $\MM{\xi}$ in a vector space $V$, subject to a  condition imposed by a \bfi{velocity map} from the vector space $V$ to the tangent space $T_{\MM{Q}}M$ of a manifold $M$ at the point $\MM{Q}$,
\begin{equation}
\mathcal{L}_\xi:\,V\times M \mapsto T_QM
\,.
\end{equation}
The velocity map $\mathcal{L}_\xi$ introduces the dynamics,
\begin{equation}
\label{q dot}
\MM{\dot{Q}}(t) = \mathcal{L}_\xi\MM{Q}(t),
\end{equation}
where $\MM{\xi}\in V$ and $\MM{\dot{Q}}\in T_QM$ is tangent to the
curve $\MM{Q}(t)$ in the manifold $M$.

Such variational principles arise in two different contexts:
\begin{enumerate}
\item The \bfi{optimal control} context, in which one seeks solutions
  for $\MM{Q}(t)$ governed by the dynamics \eqnref{q dot} that control
  the motion along a curve in an interval $0\leq t\leq T$ so as to
  minimise the cost in (\ref{lag-cost}) for a given cost functional
  $l[\MM{\xi}]$.
\item The \bfi{Hamilton's principle} context, in which stationarity
  $\delta S=0$ of the action in \eqnref{lag-cost} implies dynamical
  equations for $\MM{\xi}$ subject to the constraint imposed by the
  velocity map (\ref{q dot}).
\end{enumerate}

One approach that applies in both contexts was first introduced for
ideal fluid dynamics in Serrin \cite{Se1959} and in Seliger and
Whitham \cite{SeWh1968}. This is the \bfi{Clebsch approach} for
deriving variational principles for Eulerian fluid dynamics. A similar
approach later emerged in the work of Bloch et al. \cite{BlCrMaRa1998}
in the optimal control of rigid bodies. This approach enforces
equation \eqnref{q dot} through a Lagrange multiplier term in the
action or cost. Doing so produces dynamical equations for $\MM{Q}$ and
for the Lagrange multiplier $\MM{P}$ in terms of $\MM{{\xi}}$,
together with a formula for $\MM{{\xi}}$ given in terms of $\MM{Q}$
and $\MM{P}$. The Lagrange multiplier $\MM{P}$ is also the canonically
conjugate momentum for $\MM{Q}$ in the corresponding Hamiltonian
formulation, and the formula for $\MM{{\xi}}$ in terms of $\MM{Q}$ and
$\MM{P}$ has special significance in the Hamiltonian framework.

Such variational principles are said to be \bfi{implicit Lagrangian
  systems} and the subject has now reached a high state of
mathematical development \cite{YoMa2006}.  In this paper we take a
``bare hands'' approach to investigating this sort of problem.

Section 2 describes the conditions under which the coordinate $\MM{Q}$
and its canonical momentum $\MM{P}$ in $T^*M$ may be eliminated from
an implicit Lagrangian system in order to obtain a dynamical system
for $\MM{{\xi}}$ only. Answering this question summons a Lie algebra
structure on V. That is, $\MM{Q}$ and $\MM{P}$ may be eliminated if
and only if $\mathcal{L}_{\MM{\xi}}$ corresponds to the action of some
Lie algebra $\mathfrak{g}$ on $M$; that is
$\mathcal{L}_{\MM{\xi}}:\mathfrak{g}\times M\to TM$ for
$\MM{\xi}\in\mathfrak{g}$. In this case the dynamical system for
$\MM{{\xi}}$ is always the Euler-Poincar\'e equation for V with the
appropriate $\ad^*$ operator defined on the dual Lie algebra
$\mathfrak{g}^*$ through the natural pairing induced by taking
variational derivatives. The formula for $\MM{{\xi}}$ in terms of
$\MM{Q}$ and $\MM{P}$ is then found \emph{via} a cotangent-lift
momentum map. This key result for the continuous case is stated and
proved in Theorem \ref{reduction}.

In the case where V is the Lie algebra of vector fields
$\mathfrak{X}(\mathbb{R}^n)$, there is a choice of Lie algebra actions
on different spaces \emph{e.g.} Lie derivatives, left-action on
embedding space \emph{etc.} In each case the resulting system for
$\MM{{\xi}}$ is the Euler-Poincar\'e equation for the diffeomorphism
group. This is the EPDiff equation.

Potential energy terms may also be introduced into the action in an
implicit Lagrangian system by following standard procedure for
Hamilton's principle. This allows an extension of the Clebsch
variational principle to obtain Euler-Poincar\'e equations with
advected quantities \cite{HoMaRa1998}. Many equations of fluid
dynamics may be obtained this way. Section 3 discusses the standard
example of the incompressible Euler equations and provides references
for other examples. Again all calculations are performed explicitly.

Clebsch variational principles not only unify the subject, they also
provide a systematic framework for deriving numerical integrators. A
great deal of activity and rapid development of these variational
integrators has recently transpired. See \cite{BoMa2007} for an
up-to-date survey of the subject, a bibliography and new results from
the same viewpoint as the present paper.

Section 4 discusses the potential for constructing numerical
integration methods by discretising the Clebsch variational approach
in both space and time and deriving the resulting discrete equations
of motion. The Clebsch approach provides a method of obtaining
variational integrators \cite{LeMaOrWe2003} simply by discretising the
Clebsch variational principle in both space and time. These
integrators are symplectic and hence they fit into the backward-error
analysis framework \cite{LeRe2005}. This means that they preserve the
Hamiltonian within $\mathcal{O}(\Delta t^p)$ ($p$ is the order of the
method in time) and also preserve conservation laws associated with
symmetries (provided that the symmetries are retained by the spatial
discretisation). We show that for dynamics on finite-dimensional Lie
groups, the discrete Clebsch variational principle results in
equations where $\MM{Q}$ and $\MM{P}$ can once again be eliminated to
produce a conservative integrator for the equation for $\MM{{\xi}}$
only. We give an example of a Clebsch integrator for the rigid body
equation which after eliminating $\MM{Q}$ and $\MM{P}$ takes a
particularly elegant form. The last section of the paper closes by
discussing some possible directions for applying the Clebsch approach
and other recently developed parallel approaches for the case of
infinite-dimensional systems.

\section{Clebsch variational principles}

\begin{definition}[Velocity map]
  Suppose $V$ is a vector space and $M$ is a manifold. For each
  $\MM{Q}\in M$ and $\MM{\xi}\in V$, we define the linear
  \bfi{velocity map} $\mathcal{L}_\xi:\,V\times M \mapsto T_QM$.
\end{definition}
For a given element $\MM{\xi}$ of $V$, we use $\mathcal{L}$ to define
velocity on $T_{\MM{Q}}M$ \emph{via}
\[
\MM{\dot{Q}} = \, \mathcal{L}_{\MM{\xi}}\MM{Q}, \qquad \MM{Q}\in M.
\]
\begin{definition}[Clebsch action principle]
For a given functional $l:V\to\mathbb{R}$, the Clebsch action
principle is
\begin{equation}
\label{Clebsch action}
\delta \int_{t_1}^{t_2} l[\MM{\xi}(t)]
+ \left\langle\MM{P}(t),
\MM{\dot{Q}}(t)-\mathcal{L}_{\MM{\xi}(t)}\MM{Q}\right\rangle_{T^*M}
\diff{t}=0,
\end{equation}
where $\MM{P}$ is a Lagrange multiplier in $T^*_{\MM{Q}}M$ and
$\left\langle\cdot,\cdot\right\rangle_{T^*M}$ is the standard inner product on
$T^*M$.
\end{definition}
\begin{remark}
The solutions of
  this action principle minimise $\int_{t_1}
  ^{t_2}l(\MM{\xi})\diff{t}$ subject to the constraint that $\MM{Q}$
  is directed by the velocity map $\mathcal{L}$.
  The second term in (\ref{Clebsch action}) imposes the definition of
  velocity by a Lagrange multiplier that will turn out to be the
  conjugate momentum in the Hamiltonian formulation. This will be
  shown in theorem \ref{reduction}.
\end{remark}
To write down the
general solutions, we need to define the diamond operator.
\begin{definition}[Diamond operator]
\label{diamond}
Let $\mathcal{L}$ be a velocity map from $V$ to $M$ as defined above. The
operator $\diamond: T^*M\to V^*$ satisfies
\[
\left\langle\MM{P}\diamond\MM{Q},\MM{\xi}\right\rangle_V = 
-\left\langle \MM{P},\mathcal{L}_{\MM{\xi}}\MM{Q}\right\rangle_{T^*M}.
\]
where $\left\langle\cdot,\cdot\right\rangle_V$ is the inner product on
$V\times V^*$.
\end{definition}
\begin{remark}
For the case where the velocity map is minus the Lie derivative, the diamond
operation is the dual action of the Lie derivative.
\end{remark}
\begin{remark}
  Later we shall see that for the case where the velocity map is a Lie
  algebra action, the quantity $-\MM{P}\diamond \MM{Q}$ is a
  cotangent-lift momentum map.
\end{remark}
\begin{lemma}[Clebsch equations]$\quad$\\
The optimising solutions for the action principle (\ref{Clebsch action}) 
satisfy:
\begin{eqnarray*}
\dede{l}{\MM{\xi}} & = & -\MM{P}\diamond\MM{Q}, \\
\MM{\dot{Q}} & = & \mathcal{L}_{\MM{\xi}}\MM{Q}, \\
\MM{\dot{P}} & = & -\left(T_{\MM{Q}}\mathcal{L}_{\MM{\xi}}\right)^T\MM{P}, \\
\end{eqnarray*}
where one defines the variational derivative operations
\[
\delta l = \left\langle\dede{l}{\MM{\xi}},\delta\MM{\xi}\right\rangle_V,
\quad \mathrm{and} \quad 
\langle\MM{P},T_{\MM{Q}}\mathcal{L}_{\MM{\xi}}\delta\MM{Q}\rangle =
\langle (T_{\MM{Q}}\mathcal{L}_{\MM{\xi}})^T\MM{P},\delta\MM{Q}\rangle.
\]
\end{lemma}
\begin{proof}
\begin{eqnarray*}
0 & = & \delta \int_{t_1}^{t_2} l[\MM{\xi}] 
+ \left\langle\MM{P},\MM{\dot{Q}}-\mathcal{L}_{\MM{\xi}}\MM{Q}\right\rangle_{T^*M}
\diff{t}, \\
 & = & \int_{t_1}^{t_2}\left\langle
\dede{l}{\MM{\xi}},\delta\MM{\xi}\right\rangle_V
+ \left\langle \delta\MM{P},\MM{\dot{Q}}-\mathcal{L}_{\MM{\xi}}\MM{Q}\right\rangle_{T^*M} \\
& & + \left\langle\MM{P},\dot{\delta\MM{Q}} - (T_{\MM{Q}}\mathcal{L}_{\MM{\xi}})\cdot\delta\MM{Q}
- \mathcal{L}_{\MM{\delta\MM{\xi}}}\MM{Q}\right\rangle_{T^*M}\diff{t}, \\
& = & \left\langle\dede{l}{\MM{\xi}}+\MM{P}\diamond\MM{Q},\delta\MM{\xi}\right\rangle_V
+ \left\langle \delta\MM{P},\MM{\dot{Q}}-\mathcal{L}_{\MM{\xi}}\MM{Q}\right\rangle_{T^*M} \\
& & + \left\langle - \MM{\dot{P}}  - \left(T_{\MM{Q}}\mathcal{L}_{\MM{\xi}}\right)^T\MM{P}, 
  \delta\MM{Q}\right\rangle_{T^*M}\diff{t},
\end{eqnarray*}
and the result follows since $\delta\MM{\xi}$, $\delta\MM{P}$, and
$\delta\MM{Q}$ are arbitrary.
\qed\end{proof}
\begin{lemma}[Legendre transform]
\label{hamiltonian}
Suppose $l$ is chosen such that it is possible to solve for
$\MM{\xi}$ from $\delta{l}/\delta{\MM{\xi}}$, \emph{i.e.} there exists 
an operator $G$ such that
\[
\MM{\xi} = G\dede{l}{\MM{\xi}}.
\]
Then the Clebsch equations for $\MM{Q}$ and $\MM{P}$ are canonical
with Hamiltonian given by
\[
H = \langle \MM{P},\mathcal{L}_{G(\MM{P}\diamond\MM{Q})}\MM{Q}\rangle -
l\left[G\left(\MM{P}\diamond\MM{Q}\right)\right].
\]
\end{lemma}
\begin{proof}
The result follows directly from calculating the canonical equations for this
Hamiltonian.
\qed\end{proof}
\begin{example}[The rigid body: Bloch et al. 1998] For example,
  consider the rigid body in coordinates $(Q,\dot{Q})\in TSO(3)$. The
  angular velocity $\Omega$ of the rigid body is defined as
\[
\dot{Q} = Q\Omega
\quad\hbox{  with }\quad
\Omega^T=-\,\Omega\in \mathfrak{so}(3).
\]
The Clebsch variational principle for the rigid body is given by
$\delta S=0$ with action $S=\int L\,dt$ and implicit Lagrangian
\[
L(\Omega,P,Q)=\frac{1}{2}\Omega\cdot I\Omega + P^T\cdot (\dot{Q} - Q\Omega)
\]
for constant, symmetric $I$. The variations of $L$ are given by:
\[
\delta \int L\,dt = \int \left[ \delta\Omega\cdot (I\Omega -  Q^TP)
+ \delta P^T\cdot (\dot{Q} - Q\Omega)
- \delta Q^T\cdot (\dot{P} - P\Omega)\right]
\,dt.
\]
Thus, stationarity  of this  implicit variational principle  implies a
set of rigid body equations which  first appeared in the work of Bloch
et al.  \cite{BlCrMaRa1998} on optimal control of rigid bodies
\[
I\Omega =  -P\diamond Q = Q^TP
\,,\quad
\dot{Q} = \mathcal{L}_\Omega Q = Q\Omega
\,,\quad
\dot{P} = (T_Q\mathcal{L}_\Omega Q)^TP = P\Omega
\,.
\]

In this particular example, one can study the dynamics of $\Omega$
by calculating
\begin{eqnarray*}
\dd{}{t}I\Omega & = & \dot{Q}^TP + Q^T\dot{P}, \\
& = & (Q\Omega)^TP + Q^TP\Omega, \\
& = & \Omega^T(Q^TP) + (Q^TP)\Omega, \\
& = & \Omega^T(I\Omega) + (I\Omega)\Omega = [I\Omega,\Omega].
\end{eqnarray*}
Thus, the Clebsch equations yield the dynamics for $\Omega$
by eliminating $P$ and $Q$ \emph{via} the derivative of $\delta{l}/\delta{\xi}$.
\end{example}
We will determine when it is possible to obtain a 
dynamical system for $\MM{\xi}$ by eliminating the
Clebsch variables. First we need to make two more definitions.
\begin{definition}[Closure]
  Let $V$ be a vector space with $\mathcal{L}$ being a velocity map
  $V\times M\to TM$.  The velocity map $\mathcal{L}$ is said to be
  \bfi{closed} if, for every pair of vectors $\MM{u},\MM{v}\in V$,
  there exists a third vector $\MM{w}\in V$ such that
\[
\mathcal{L}_{\MM{w}}\MM{Q} = ((T_{\MM{Q}}\mathcal{L}_{\MM{v}})\mathcal{L}_{\MM{u}}
-(T_{\MM{Q}}\mathcal{L}_{\MM{u}})\mathcal{L}_{\MM{v}})\MM{Q},
\]
where $T_{\MM{Q}}\mathcal{L}_{\MM{u}}$ is the tangent of $\mathcal{L}_{\MM{u}}$ evaluated
at $\MM{Q}$.  In that case, the velocity map induces an algebra
structure on $V$ with bracket
\[
\mathcal{L}_{[\MM{u},\MM{v}]} = (T_{\MM{Q}}\mathcal{L}_{\MM{v}})\mathcal{L}_{\MM{u}}
-(T_{\MM{Q}}\mathcal{L}_{\MM{u}})\mathcal{L}_{\MM{v}},
\]
which satisfies the Jacobi identity.
\end{definition}
\begin{definition}[$\ad$- and $\ad^*$-operators]
Let $\mathcal{L}$ define an velocity map from $V$ on $M$ with induced bracket
$[\cdot,\cdot]:V\times V\to V$. We define the $\ad$-operator by
\[
\ad_{\MM{u}}\MM{v} = -[\MM{u},\MM{v}].
\]
For $\MM{m}$ in $V^*$ we define the $\ad^*$-operator by 
\[
\left\langle \ad^*_{\MM{u}}\MM{m},\MM{v}\right\rangle_V = 
\left\langle \MM{m},\ad_{\MM{u}}\MM{v}\right\rangle_V. 
\]
\end{definition}
\begin{theorem}[Elimination requires closure of the velocity map]
$\quad$\\
\label{reduction}
Let $V$ be a vector space, and let $\mathcal{L}$ be a velocity map
from $V\times M\to TM$. Then the cotangent variables
$(\MM{P},\MM{Q})\in T^*M$ may be eliminated from the Clebsch equation
if and only if the image of $\mathcal{L}$ in $TM$ is closed under the Lie
bracket.
\smallskip
\noindent Furthermore, when this closure condition holds:
\begin{enumerate}
\item $\mathcal{L}$ induces a Lie algebra structure on $V$,
\item $-\MM{P}\diamond\MM{Q}$ is a cotangent-lift 
momentum map, and
\item $\MM{\xi}$ satisfies the 
Euler-Poincar\'e equation:
\[
\dd{}{t}\dede{l}{\MM{\xi}} + \ad_{\MM{\xi}}^*\dede{l}{\MM{\xi}}=0.
\]
\end{enumerate}
\end{theorem}
\begin{proof}
First suppose that $\mathcal{L}$ is closed. For any vector $\MM{w}\in V$,
\begin{eqnarray*}
\dd{}{t}\left\langle\dede{l}{\MM{\xi}}(\MM{\xi}),\MM{w}\right\rangle_V 
& = & -\,\dd{}{t}\langle\MM{P}\diamond\MM{Q},\MM{w}\rangle_V \\
\hbox{[Definition of $\diamond$]}  & = & \dd{}{t}\langle\MM{P},\mathcal{L}_{\MM{w}}\MM{Q}\rangle_{T^*M} \\
\hbox{[$\MM{P}$ and 
$\MM{Q}$ equations]} & = & \langle -(T_{\MM{Q}}\mathcal{L}_{\MM{\xi}})^T\MM{P},\mathcal{L}_{\MM{w}}\MM{Q}\rangle_{T^*M} 
 + \langle \MM{P},(T_{\MM{Q}}\mathcal{L}_{\MM{w}})\mathcal{L}_{\MM{\xi}}\MM{Q}\rangle_{T^*M} \\& = & \langle \MM{P}, \left(-(T_{\MM{Q}}\mathcal{L}_{\MM{\xi}})\mathcal{L}_{\MM{w}} + (T_{\MM{Q}}\mathcal{L}_{\MM{w}})\mathcal{L}_{\MM{\xi}}\right)\MM{Q}\rangle_{T^*M} \\
\hbox{[closure]} & = & \langle \MM{P},-\mathcal{L}_{[\MM{\xi},\MM{w}]}\MM{Q}\rangle_{T^*M} \\
\hbox{[Definition of $\diamond$]}
& = & \langle \MM{P}\diamond\MM{Q}, [\MM{\xi},\MM{w}] \rangle_V \\
& = & -\left\langle \dede{l}{\MM{\xi}}(\MM{\xi})
, \ad_{\MM{\xi}}\MM{w} \right\rangle_V, \\
& = & -\left\langle \ad^*_{\MM{\xi}}\dede{l}{\MM{\xi}}(\MM{\xi})
, \MM{w} \right\rangle_V, \\
\end{eqnarray*}
and hence
\[
\dd{}{t}\dede{l}{\MM{\xi}} + \ad^*_{\MM{\xi}}\dede{l}{\MM{\xi}} = 0,
\]
which is a closed system for $\MM{\xi}$. Hence, closure allows elimination
of $\MM{P}$ and $\MM{Q}$.

Conversely, suppose that $\mathcal{L}$ is not closed. Then there
exist $\MM{u},\MM{v}\in V$ such that
\begin{equation}
  (T_{\MM{Q}}\mathcal{L}_{\MM{v}})\mathcal{L}_{\MM{u}}-(T_{\MM{Q}}\mathcal{L}_{\MM{u}})\mathcal{L}_{\MM{v}} \neq
  \mathcal{L}_{\MM{w}}\,,
\label{assume}
\end{equation}
for any $\MM{w}\in V$. Now assume, aiming for a contradiction, that
$\MM{P}$ and $\MM{Q}$ may be eliminated from the equations. For any
$\MM{Q}$, we may find a $\MM{P}$ such that
\[
\dede{l}{\MM{\xi}}(\MM{u}) 
:= \dede{l}{\MM{\xi}}\Big|_{\MM{\xi}=\MM{u}}
 = \MM{P}\diamond\MM{Q}.
\]
Then
\begin{eqnarray}
\left\langle\dd{}{t}\dede{l}{\MM{u}}(\MM{u}),\MM{v}\right\rangle_V 
& = & \langle \MM{P}, \left(-(T_{\MM{Q}}\mathcal{L}_{\MM{u}})\mathcal{L}_{\MM{v}} 
+ (T_{\MM{Q}}\mathcal{L}_{\MM{v}})\mathcal{L}_{\MM{u}}\right)\MM{Q}\rangle. 
\label{mv-pairing}
\end{eqnarray}
In order for this to be consistent, the left-hand side must be linear
in $\MM{\xi}$ and $\delta l/\delta\MM{\xi}$ (treating $\delta l/\delta
\MM{\xi}$ as a separate variable), and hence we may write
\[
\dd{}{t}\dede{l}{\MM{\xi}}(\MM{u}) =
\hat{\mathcal{L}}_{\MM{u}}\dede{l}{\MM{\xi}}(\MM{u}),
\]
for some linear operator $\hat{\mathcal{L}}_{\MM{u}}$. We pair this
with a vector, then use the Clebsch variational equation and the
definition of diamond to write
\begin{eqnarray*}
\left\langle\dd{}{t}\dede{l}{\MM{\xi}}(\MM{u}),\MM{v}\right\rangle_V
 & = & \left\langle\hat{\mathcal{L}}_{\MM{u}}
\dede{l}{\MM{\xi}}(\MM{u}),\MM{v}\right\rangle_V \\
& = & \left\langle
\dede{l}{\MM{\xi}}(\MM{u}),\hat{\mathcal{L}}_{\MM{u}}^*\MM{v}\right\rangle_V \\
& = & -\left\langle
\MM{P}\diamond\MM{Q},\hat{\mathcal{L}}_{\MM{u}}^*\MM{v}\right\rangle_V \\
& = & \left\langle
\MM{P},\mathcal{L}_{\hat{\mathcal{L}}_{\MM{u}}^*\MM{v}}\MM{Q}\right\rangle_V.
\end{eqnarray*}
Consequently, comparing with equation (\ref{mv-pairing}) yields
\[
\mathcal{L}_{\hat{\mathcal{L}}_{\MM{u}}^*\MM{v}}\MM{Q} =
[ (T_{\MM{Q}}\mathcal{L}_{\MM{u}})\mathcal{L}_{\MM{v}}-(T_{\MM{Q}}\mathcal{L}_{\MM{v}})\mathcal{L}_{\MM{u}}]\MM{Q}
\]
for any $\MM{Q}$, which contradicts the assumption
(\ref{assume}). Hence, $\MM{P}$ and $\MM{Q}$ may be eliminated from
the Clebsch equation if and only if the image of $\mathcal{L}$ in $TM$
is closed under the Lie bracket.

The fact that $-\MM{P}\diamond\MM{Q}$ is a cotangent-lift momentum
map comes straight from the definition: if $\mathcal{L}$ is closed, then $V$ is
a Lie algebra with bracket induced by $\mathcal{L}$, and the map $\mathcal{L}$ is a Lie
algebra action on $M$. The Hamiltonian for the corresponding
cotangent-lifted action is
\[
h_{\MM{\xi}} = \langle \MM{P}, \mathcal{L}_{\MM{u}}\MM{Q}\rangle_{T^*M},
\]
and the momentum map $\MM{J}$ for this action is defined by the relation
\[
\langle\MM{J},\MM{\xi}\rangle_V = \langle \MM{P},\mathcal{L}_{\MM{u}}\MM{Q}\rangle_{T^*M}
\]
which matches Definition \ref{diamond}, and hence
\[
\MM{J} = -\MM{P}\diamond\MM{Q}.
\]
\qed\end{proof}
\begin{remark}
As a result of Theorem \ref{reduction}, solutions of the Clebsch equations
may be composed by $\diamond$ into solutions of the following equation 
\[
\dd{}{t}\dede{l}{\MM{\xi}} + \ad_{\MM{\xi}}^*\dede{l}{\MM{\xi}} = 0,
\]
whenever $\mathcal{L}$ is closed. This equation is called the Euler-Poincar\'e
equation and it describes geodesic motion on the space $V$ whenever
$l$ is a quadratic functional (metric).
\end{remark}
To extend the scope of this framework it is useful to introduce potential
energy terms to $l$ in the action principle, \emph{i.e.} let $l$ be a function
of the elements of $M$ as well as $V$.
\begin{definition}[Clebsch action principle with potential energy]
$\quad$\\
For a given functional $l:V\times M\to\mathbb{R}$, the Clebsch action
principle is
\begin{equation}
\label{Clebsch action with potential}
\delta \int_{t_1}^{t_2} l[\MM{\xi}(t),\MM{Q}(t)] +
\left\langle\MM{P}(t),\MM{\dot{Q}}(t)-\mathcal{L}_{\MM{\xi}(t)}\MM{Q}\right\rangle_{T^*M}
\diff{t}=0,
\end{equation}
where $\MM{P}$ is a Lagrange multiplier in $T^*_{\MM{Q}}M$ and
$\left\langle\cdot,\cdot\right\rangle_{T^*M}$ is the usual inner product on
$T^*M$.
\end{definition}
\begin{lemma}[Clebsch equations with potential energy]
The optimising solutions for the action principle (\ref{Clebsch action}) are:
\begin{eqnarray*}
\dede{l}{\MM{\xi}} & = & -\MM{P}\diamond\MM{Q}, \\ \MM{\dot{Q}} & = &
\mathcal{L}_{\MM{\xi}}\MM{Q}, \\ \dot{\MM{P}} & = &
-\left(T_{\MM{Q}}\mathcal{L}_{\MM{\xi}}\right)^T\MM{P} + \dede{l}{\MM{Q}}. \\
\end{eqnarray*}
\end{lemma}
\begin{proof}
\begin{eqnarray*}
0 & = & \delta \int_{t_1}^{t_2} l[\MM{\xi},\MM{Q}] 
+ \left\langle\MM{P},\MM{\dot{Q}}-\mathcal{L}_{\MM{\xi}}\MM{Q}\right\rangle_{T^*M}
\diff{t}, \\
 & = & \int_{t_1}^{t_2}\left\langle
\dede{l}{\MM{\xi}},\delta\MM{\xi}\right\rangle_V
+ \left\langle \dede{l}{\MM{Q}},\delta\MM{Q}\right\rangle
+ \left\langle \delta\MM{P},\MM{\dot{Q}}-\mathcal{L}_{\MM{\xi}}\MM{Q}\right\rangle_{T^*M} \\
& & + \left\langle\MM{P},\dot{\delta\MM{Q}} - (T_{\MM{Q}}\mathcal{L}_{\MM{\xi}})\cdot\delta\MM{Q}
- \mathcal{L}_{\MM{\delta\MM{\xi}}}\MM{Q}
\right\rangle_{T^*M}\diff{t}, \\
& = & \left\langle\dede{l}{\MM{\xi}}+\MM{P}\diamond\MM{Q},\delta\MM{\xi}\right\rangle_V
+ \left\langle \delta\MM{P},\MM{\dot{Q}}-\mathcal{L}_{\MM{\xi}}\MM{Q}\right\rangle_{T^*M} \\
& & + \left\langle -\dot{\MM{P}} -
\left(T_{\MM{Q}}\mathcal{L}_{\MM{\xi}}\right)^T\MM{P} + \dede{l}{\MM{Q}},
\delta\MM{Q}\right\rangle_{T^*M}\diff{t},
\end{eqnarray*}
and the result follows since $\delta\MM{\xi}$, $\delta\MM{P}$, and
$\delta\MM{Q}$ are arbitrary.
\qed\end{proof}
\begin{theorem}[Elimination theorem with potential energy]
Let $V$ be a vector space, and let $\mathcal{L}$ be a velocity map from $V\to TM$. The cotangent variables $\MM{P}$ and $\MM{Q}$ may be eliminated from the Clebsch
equation with potential energy, if and only if the image of 
$\mathcal{L}$ in $TM$ is closed under the Lie bracket.

Furthermore, when the condition holds, $\MM{\xi}$ and $\MM{Q}$ satisfy
the Euler-Poincar\'e equation with advected quantities:
\[
\dd{}{t}\dede{l}{\MM{\xi}} + \ad_{\MM{\xi}}^*\dede{l}{\MM{\xi}}=
-\dede{l}{\MM{Q}}\diamond\MM{Q}, 
\quad \MM{Q}_t = \mathcal{L}_{\MM{\xi}}\MM{Q}.
\]
\end{theorem}
\begin{proof}
The proof follows the proof of Theorem \ref{reduction}, using
\begin{eqnarray*}
\dd{}{t}\Bigg\langle \dede{l}{\MM{\xi}},\MM{w}\Bigg\rangle_V & = & 
-\,\dd{}{t}\Bigg\langle
\MM{P}\diamond\MM{Q} 
,\MM{w}
\Bigg\rangle_V,
\\ 
& = & \dd{}{t}\Bigg\langle \MM{P},\mathcal{L}_{\MM{w}}\MM{Q}\Bigg\rangle_V  \\
& = & \left\langle -(T_{\MM{Q}}\mathcal{L}_{\MM{\xi}})^T\MM{P} + \dede{l}{\MM{Q}}
,\mathcal{L}_{\MM{w}}\MM{Q}\right\rangle_{T^*M} 
 + \Bigg\langle \MM{P},(T_{\MM{Q}}\mathcal{L}_{\MM{w}})\mathcal{L}_{\MM{\xi}}\MM{Q}\Bigg\rangle_{T^*M}, \\
& = & \left\langle -\ad_{\MM{\xi}}^*\dede{l}{\MM{Q}} -
\dede{l}{\MM{Q}}\diamond\MM{Q},\MM{w}
\right\rangle.
\end{eqnarray*}
\qed\end{proof}

\section{Examples}
\subsection{Singular solutions of EPDiff}
When $V$ is the space of vector fields $\mathfrak{X}(\mathbb{R}^n)$,
under the conditions of Theorem \ref{reduction} the Clebsch implicit
variational principle yields the Euler-Poincar\'e equation for
diffeomorphisms (EPDiff):
\[
\dd{}{t}\dede{l}{\MM{u}} + \nabla\cdot(\MM{u}\dede{l}{\MM{u}})
+ (\nabla\MM{u})^T\cdot\dede{l}{\MM{u}} = 0.
\]

One possible way to form a Clebsch principle for EPDiff is to consider the 
left-action of vector fields on the space of embeddings
$M=\Emb(S,\mathbb{R}^n)$
for some manifold $S$ (such as the circle, or the sphere).
For an embedding
\[
\MM{Q}:S\to\mathbb{R}^n,
\]
the velocity map is defined by composition of functions
\[
\mathcal{L}_{\MM{u}}\MM{Q}(s) = \MM{u}(\MM{Q}(s)), \quad s\in S.
\]
The Clebsch principle is then
\[
\delta\int_{t_1}^{t_2} l[\MM{u}] + \langle
\MM{P},\MM{\dot{Q}}-\MM{u}(\MM{Q}) \rangle_{T^*M}\diff{t}=0,
\]
where the inner product in the second term is defined as
\[
\int_S \MM{P}(s,t)\cdot\left
(\MM{\dot{Q}}(s,t)-\MM{u}(\MM{Q}(s,t))\right)\diff s.
\]
The diamond operator $(\diamond)$ is thus defined in this case by
\begin{eqnarray*}
\langle\left(\MM{P}\diamond\MM{Q}\right),\MM{u}\rangle_V
 &=& \langle \MM{P}, \mathcal{L}_{\MM{u}}\MM{Q} \rangle_{T^*M} \\
 &=& \int_S \MM{P}(s)\cdot\MM{u}(\MM{Q}(s))\diff s, \\
 & = & \int_S \MM{P}(s)
\int_M\delta(\MM{x}-\MM{Q}(s,t))\MM{u}(\MM{x})\diff{V}(\MM{x})\diff s \\
 & = & \int_M\left(\int_S \MM{P}(s)
\delta(\MM{x}-\MM{Q}(s,t))
\diff s\right)\cdot \MM{u}(\MM{x})\diff{V}(\MM{x})
.
\end{eqnarray*}
Consequently, one finds
\[
\MM{P}\diamond\MM{Q} = \int_S \MM{P}(s)
\delta(\MM{x}-\MM{Q}(s,t))
\diff s
\,,
\]
which is the singular solution ansatz of \cite{HoMa2004}.

One then calculates the Clebsch equations as 
\begin{eqnarray}
\dede{l}{\MM{u}} &=& \MM{P}\diamond\MM{Q}
=\int_S  \MM{P}(s,t)\delta(\MM{x}-\MM{Q}(s,t))
\diff s, \label{singular}\\
\MM{\dot{Q}} & = & \mathcal{L}_{\MM{u}}\MM{Q} = \MM{u}(\MM{Q}), \\
\dot{\MM{P}} & = & -(T_{\MM{Q}}\mathcal{L}_{\MM{u}})^T\MM{P} = 
-\sum_kP_k\pp{u^k}{\MM{Q}}=-(\nabla\MM{u}(\MM{Q}))^T\cdot\MM{P}, 
\end{eqnarray}
which we know to be canonically Hamiltonian from Lemma
\ref{hamiltonian}.  Furthermore, we know that
$\delta{l}/\delta{\MM{u}}$ satisfies the EPDiff equation from Theorem
\ref{reduction}. To verify this statement, take the inner product with a test
function and differentiate in time, as follows.
\begin{eqnarray*}
\dd{}{t}\left\langle \dede{l}{\MM{u}}, \MM{w}\right\rangle_{\mathfrak{X}(M)} & = & 
\dd{}{t} \langle \MM{P}\diamond \MM{Q}, \MM{w} \rangle_{\mathfrak{X}(M)}, \\
 & = & \dd{}{t} \int_M
\left( \int_S \MM{P}(s,t)\delta(\MM{x}-\MM{Q}(s,t))\diff{s}\right)\cdot \MM{w}(\MM{x})\diff{\Vol}
(\MM{x}), \\
 & = & \dd{}{t} \int_S \MM{P}\cdot \MM{w}(\MM{Q}) \diff{s}, \\
 & = & \int_S(\dot{\MM{P}}\cdot \MM{w}(\MM{Q}) + \MM{P} \cdot \pp{\MM{w}}{\MM{Q}}\cdot\MM{\dot{Q}})\diff{s}, \\
 & = & \int_S -\left((\nabla\MM{u})^T\cdot\MM{P}\right)\cdot \MM{w}(\MM{Q}) + \MM{P}\cdot\nabla \MM{w}(\MM{Q}) \cdot \MM{u}(\MM{Q}) 
\diff{s}, \\
& = & \int_M \Bigg(
-\int_S \MM{P}\delta(\MM{x}-\MM{Q}(s))\diff s\cdot \nabla \MM{u}(\MM{x})\cdot \MM{w}(\MM{x}), \\
& & \qquad 
+ \int_S \MM{P}\delta(\MM{x}-\MM{Q}(s))\diff{s} \cdot \nabla \MM{w}(\MM{x}) \cdot \MM{u}(\MM{x})\Bigg)
\diff{\Vol}(\MM{x}), \\
& = & -\left\langle (\nabla\MM{u})^T\cdot\dede{l}{\MM{u}},
\MM{w}\right\rangle_{\mathfrak{X}(M)} +
\left\langle \dede{l}{\MM{u}}, \nabla \MM{w}\cdot \MM{u}\right\rangle_{\mathfrak{X}(M)}.
\end{eqnarray*}
This shows that $\delta{l}/\delta\MM{u}$ satisfies the weak form of EPDiff:
\[
\dd{}{t}\left\langle \dede{l}{\MM{u}}, \MM{w}\right\rangle_{\mathfrak{X}(M)}
+ \left\langle (\nabla \MM{u})^T\cdot\dede{l}{\MM{u}}
, \MM{w}\right\rangle_{\mathfrak{X}(M)}
+ \left\langle \dede{l}{\MM{u}}, \nabla \MM{w} \cdot \MM{u}\right\rangle_{\mathfrak{X}(M)}
= 0.
\]

\begin{remark}
  In consonance with Theorem \ref{reduction} and as discussed in
  \cite{HoMa2004} the singular solution ansatz for EPDiff given by
  \eqnref{singular} is an equivariant momentum map arising from the
  cotangent lift of the action of vector fields corresponding to
  composition the left.
\end{remark}

\subsection{Euler equations}
The Lagrangian for the incompressible Euler equations is 
 \[
l[\MM{u},D] = \int_M \frac{D}{2}|\MM{u}|^2 + p(1-D)\diff{\Vol}(\MM{x}),
\]
where $D$ is the density, and $p$ is a Lagrange multiplier introduced
to enforce the constraint that the fluid is incompressible. Not
unexpectedly, the quantity $p$ turns out to be the pressure.

There are several ways to write down a Clebsch variational principle
for the incompressible Euler equations. However in order to obtain a
set of variables which include all possible solutions one needs to
include at least $d$ scalar Lagrange multipliers where $d$ is the
spatial dimension of $M$. If, for example, we only use one Lagrange
multiplier then only the vorticity-free solutions are obtained, as
first noticed by Lin (see \cite{CeMa1987} for discussion and
references).

As described in \cite{HoKu1983}, one possible way to construct such a
Clebsch variational principle is to use the action of the
diffeomorphisms on the Lagrangian labels, which satisfy
\[
\ell^A_t 
= \mathcal{L}_{\MM{u}}\ell^A 
=  -\,\MM{u}\cdot\nabla \ell^A
\,, 
\quad \ell^A(\MM{x},0)= x^A, \quad A=1,\ldots,d.
\]
The Clebsch variational principle is then
\begin{eqnarray*}
\delta \int_{t_1}^{t_2}
\int_M \frac{D}{2}|\MM{u}|^2 - p(D-1)\diff{\Vol}(\MM{x}) 
& & \\
+
\int_M
\MM{P}\cdot\left(
\MM{\ell}_t+\MM{u}\cdot\nabla\MM{\ell}\right)
\diff{\Vol}(\MM{x}) 
& & \\
+
\int_M
\phi\cdot\left(
D_t+\nabla\cdot(\MM{u}D)
\right)\diff{\Vol}(\MM{x})
\diff{t} 
\ & = &\ 0\,.
\end{eqnarray*}
Here we have also introduced an additional Lagrange multiplier $\phi$ to
enforce the dynamics of the density $D$ which means we have $d+1$
scalar Lagrange multipliers. The system could be reduced to $d$
multipliers but this would unnecessarily complicate the exposition.

The Clebsch equations are:
\begin{eqnarray*}
\pp{\MM{\ell}}{t} + \MM{u}\cdot\nabla\MM{\ell} & = & 0, \\
\pp{D}{t} + \nabla\cdot(\MM{u}D) & = & 0, \\
\pp{\MM{P}}{t} + \nabla\cdot(\MM{u}\MM{P}) & = & 0, \\
\pp{\phi}{t} + \MM{u}\cdot\nabla\phi & = & \dede{l}{D} = \frac{1}{2}|\MM{u}|^2
- p, \\
\dede{l}{\MM{u}} &=& D\MM{u} 
= -(\nabla\MM{\ell})^T\cdot\MM{P} - D\nabla\phi, \\
\end{eqnarray*}
together with the constraint $D=1$. 
After elimination one obtains
\begin{eqnarray*}
\pp{(\MM{u}D)}{t} 
+ \MM{u}\cdot\nabla(\MM{u}D) 
+ (\nabla\MM{u})^T\cdot\MM{u}D
&=& 
D\nabla\frac{1}{2}|\MM{u}|^2 
- \nabla p, 
\\
\pp{D}{t}
 + 
\nabla\cdot(\MM{u}D) 
&=&
 0
 \,,\quad 
 D=1
\end{eqnarray*}
which becomes
\[
\pp{\MM{u}}{t} 
+ \MM{u}\cdot\nabla\MM{u} 
= -\nabla p, \quad \nabla\cdot\MM{u}
=0,
\]
after substituting $D=1$.

The Clebsch variational principle encompasses essentially all fluid
theories.  See \cite{Ho2001} for references and the example of a
complex fluid.

\section{Clebsch integrators}

In this section we discuss the potential for constructing numerical
integration methods by discretising the Clebsch variational principle
in both space and time and thereby deriving the resulting discrete
equations of motion. Any numerical method obtained this way will
automatically be a variational integrator, and hence will inherit the
accompanying conservative properties such as exact preservation of
momenta where the discretisation preserves the symmetry of the
continuous system, and the long-time approximate conservation of
energy \emph{via} backward-error analysis. See \cite{LeMaOrWe2003} for
a review of variational integrator methods.

In the finite dimensional case one simply needs to follow the
variational integrator programme by finding a map which approximates
the time-derivative, and then substituting it into the Clebsch
variational principle. In the following sections we show how to do
this for the case of finite dimensional Lie groups such as $SO(3)$.
\begin{definition}
For a manifold $M$, define the \emph{derivative map} $\phi_{\Delta t}$
by
\[
\phi_{\Delta t}: M\times M\to TM, \qquad
(\MM{Q}^{n+1},\MM{Q}^n)\mapsto T_{\MM{Q}^{n+1}}M,
\]
such that $\phi_{\Delta t}(\MM{Q}^{n+1},\MM{Q}^n)$ is an approximation
to $\MM{\dot{Q}}$ at $\MM{Q}^n$ \cite{LeMaOrWe2003}.
\end{definition}
\begin{definition}[First-order discrete Clebsch action principle]
$\quad$\\
\label{fo discrete action}
For a given functional $l(\MM{\xi})$, the first-order discrete
Clebsch action principle is
\[
\delta F = 
\delta \sum_{n=0}^{N-1}
\left(
l(\MM{\xi}^n) + \langle\MM{P}^n,\phi_{\Delta t}(\MM{Q}^{n+1},\MM{Q}^n)
-L_{\MM{\xi}^n}\MM{Q}^n\rangle
\right)=0.
\]
\end{definition}
\begin{lemma}
The first-order discrete action principle given in Definition \ref{fo
discrete action} is optimised by $\MM{P}$, $\MM{Q}$ and $\MM{\xi}$ 
satisfying the \bfi{Discrete Clebsch equations}
\begin{eqnarray*} 
\dede{l}{\MM{\xi}}(\MM{\xi})^n + \MM{P}^n\diamond\MM{Q}^n & = & 0, \\
\phi_{\Delta t}(\MM{Q}^{n+1},\MM{Q}^n) & = & L_{\MM{\xi}^n}\MM{Q}^n, \\
\left(D_1\phi_{\Delta t}(\MM{Q}^n,\MM{Q}^{n-1})\right)^T\MM{P}^{n-1}
+\left(D_2\phi_{\Delta t}(\MM{Q}^{n+1},\MM{Q}^n)\right)^T\MM{P}^n
& = & T_{\MM{Q}^n}\left(\mathcal{L}_{\MM{\xi}^n}\MM{Q}\right)^T\MM{P}^n.
\end{eqnarray*}
\end{lemma}
\begin{proof}
\begin{eqnarray*}
\delta F & = & \sum_{n=0}^{N-1}
\Bigg(
\left\langle \dede{l}{\MM{\xi}}(\MM{\xi})^n, \delta\MM{\xi}^n
\right\rangle
- \langle \MM{P}^n, L_{\delta \MM{\xi}^n}\MM{Q}^n\rangle \\
& & \quad + \langle \delta\MM{P}^n, \phi_{\Delta t}(\MM{Q}^{n+1},\MM{Q}^n)
- L_{\MM{\xi}^n}\MM{Q}^n\rangle \\
& & \quad + \Big\langle \MM{P}^n, 
D_1\phi_{\Delta t}(\MM{Q}^{n+1},\MM{Q}^n)\delta\MM{Q}^{n+1}
+ D_2\phi_{\Delta t}(\MM{Q}^{n+1},\MM{Q}^n)\delta\MM{Q}^n \\
& & \qquad - T_{\MM{Q}^n}\left(\mathcal{L}_{\MM{\xi}^n}\MM{Q}\right)
\delta\MM{Q}^n\Big\rangle\Bigg),
\end{eqnarray*}
and the result follows after renumbering indices.
\qed\end{proof}
\begin{remark}
When $M$ is a vector space, one may choose 
\[
\phi_{\Delta t}(\MM{Q}^{n+1},\MM{Q}^n) = \frac{\MM{Q}^{n+1}-\MM{Q}^n}
{\Delta t},
\]
and the equations become
\begin{eqnarray*}
\dede{l}{\MM{\xi}}(\MM{\xi}^n) & = & -\MM{P}^n\diamond\MM{Q}^n, \\
\MM{Q}^{n+1} & = & \MM{Q}^n + \Delta tL_{\MM{\xi}^n}\MM{Q}^n, \\
\MM{P}^{n+1} & = & \MM{P}^n - \Delta t \left(T_{\MM{Q}^{n+1}}
\mathcal{L}_{\MM{\xi}^{n+1}}\MM{Q}\right)^T\MM{P}^{n+1}.
\end{eqnarray*}
\end{remark}
\begin{remark}
These equations give the first-order symplectic method, known in
\cite{LeRe2005} as symplectic Euler-A, for the Hamiltonian system
given in Lemma \ref{hamiltonian}. The adjoint method to this, known in
\cite{LeRe2005} as symplectic Euler-B, is obtained from the following
discrete variational principle:
\[
\delta \sum_{n=1}^N
\left(l(\MM{\xi}^n) + \langle \MM{P}^n,
-\phi_{\Delta t}(\MM{Q}^{n-1},\MM{Q}^n)
-\mathcal{L}_{\MM{\xi}^n}\MM{Q}^n
\rangle\right)=0.
\]
\end{remark}
\begin{remark}
  Higher-order schemes may be obtained by replacing the first-order
  discretisation of the $\MM{Q}$-equation enforced by the Lagrange
  multiplier $\MM{P}$ by Munthe-Kaas methods \cite{MK1998}
  (Runge-Kutta methods on Lie groups).
\end{remark}

\subsection{First-order integrators for matrix groups}
In this section we show how to construct a derivative map for the case
of matrix groups. For this entire section:

\begin{enumerate}
\item $\mathrm{Q}$ is a $d$-dimensional complex matrix.
\item The velocity map acts by
matrix multiplication by $\mathrm{X}$ from the right 
\[
\mathcal{L}_{\mathrm{X}}\mathrm{Q} = \mathrm{Q}\mathrm{X}.
\]
\item The inner product on $TM$ is defined by the matrix trace using
  the transpose-conjugate operation $(\dagger)$
\[
\langle\mathrm{P},\dot{\mathrm{Q}}\rangle = \Tr\left(
\mathrm{P}\mathrm{Q}^\dagger\right) = \sum_{kl}P_{kl}\bar{Q}_{kl},
\]
where the overbar indicates complex conjugation.
\item The diamond operation is given by
$\mathrm{P}\diamond\mathrm{Q}=-\,\mathrm{Q}^\dagger\mathrm{P}$ since
\begin{eqnarray*}
\langle\mathrm{P}\diamond\mathrm{Q},\mathrm{X}\rangle & = &
-\langle \mathrm{P},\mathcal{L}_{\mathrm{X}}\mathrm{Q}\rangle  \\
&=& -\Tr\left(\mathrm{P}(\mathrm{Q}\mathrm{X})^\dagger\right) \\
&=& -\Tr\left(\mathrm{Q}^\dagger\mathrm{P}\mathrm{X}^\dagger\right) \\
& = & -\langle \mathrm{Q}^\dagger\mathrm{P},\mathrm{X}\rangle.
\end{eqnarray*}
\end{enumerate} 
We first define the exponential map
when the velocity map $\mathcal{L}_{\mathrm{X}}$ is defined
by right matrix multiplication.
\begin{definition}[exponential map]
The \emph{exponential map} corresponding to the velocity map 
$\mathcal{L}_{\mathrm{X}}$ is the solution of the equation
\[
\dd{}{t}\exp(t\mathrm{X}) = \mathcal{L}_{\mathrm{X}}\exp(t\mathrm{X}):=
\exp(t\mathrm{X})\mathrm{X}, \quad
\exp(0) = \Id.
\]
\end{definition}
We also define the logarithm map:
\begin{definition}[logarithm map]
Let $\mathcal{L}$ be a velocity map. The logarithm map is defined by
\[
\log(\exp(\mathrm{X})) = \mathrm{X},
\]
for any vector element $\mathrm{X}$ of $V$.
\end{definition}
For a pair of group elements $\mathrm{Q}_1$, $\mathrm{Q}_2$, we seek a vector
$\mathrm{X}$ such that
\[
\mathrm{Q}_2 = \mathrm{Q}_1\exp(\Delta t\mathrm{X}).
\]
This allows one to use $\mathcal{L}_{\mathrm{X}}\mathrm{Q}_2$ as an
approximation for the time derivative at $\mathrm{Q}_2$, thereby motivating the following definition.
\begin{definition}[Discrete approximation of time derivative]
\label{Discretetimederiv}
$\quad$\\
For two group elements $\mathrm{Q}_1$, $\mathrm{Q}_2$ an approximation to the
time derivative at $\mathrm{Q}_1$ of a solution which passes between
$\mathrm{Q}_1$ and $\mathrm{Q}_2$ in time $\Delta t$ is 
\[
\phi_{\Delta t}(\mathrm{Q}_2,\mathrm{Q}_1) = \mathrm{Q}_1\frac{1}{\Delta t}
\logapp\left(\mathrm{Q}_1^{-1}\mathrm{Q}_2\right),
\]
where $\logapp:\Omega\to V$ is given by
\[
\logapp(A) = \log(A) + \mathcal{O}((A-I)^p),
\]
for some positive integer $p$.
\end{definition}
Let us now construct the Clebsch integrator.
\begin{definition}[First-order discretisation]
\label{euler}
We replace the time integral in the Clebsch variational principle
by a sum and substitute our discrete approximation of the time derivative 
to get
\[
\delta \sum_{n=1}^N
\left(l\left(\mathrm{X}^{n-1}\right) + \left\langle \mathrm{P}^{n-1},
\frac{1}{\Delta t}
\mathrm{Q}^{n-1}\logapp
\left(\left(\mathrm{Q}^{n-1}\right)^{-1}\mathrm{Q}^n\right)
-\mathrm{Q}^{n-1}\mathrm{X}^{n-1}
\right\rangle\right)=0.
\]
Note that the inner product is taken at the point $\mathrm{Q}^{n-1}$.
\end{definition}
\begin{theorem}
\label{discrete Clebsch}
The Clebsch variational principle given in Definition \ref{euler}
leads to the following first-order symplectic Euler discretisation:
\begin{eqnarray}
\label{discrete dede l}
\dede{l}{\mathrm{X}}(\mathrm{X}^{n-1}) & = & (\mathrm{Q}^{n-1})^\dagger\mathrm{P}^{n-1}, \\
\label{discrete Q}
\mathrm{Q}^n & = & \mathrm{Q}^{n-1}\expapp(\Delta t{\mathrm{X}^{n-1}})
, \\
\mathrm{P}^n
& = &
\left(\left(\mathrm{Q}^n\right)^\dagger\right)^{-1} 
\left(\left(T_{\expapp(\Delta t\mathrm{X}^n)}\logapp\right)^\dagger\right)^{-1}
\left(\expapp(\Delta t\mathrm{X}^{n-1})\right)^\dagger \nonumber \\
& & \qquad \left(T_{\expapp(\Delta t\mathrm{X}^{n-1})}\logapp\right)^\dagger
\left(\mathrm{Q}^{n-1}\right)^\dagger\mathrm{P}^{n-1}\left(\expapp(\Delta t\mathrm{X}^n)\right)^{-1},
\end{eqnarray}
where $\expapp$ is the inverse of the $\logapp$ operation.
\end{theorem}
\begin{remark}
The discrete Clebsch momentum map takes the expected form.
\end{remark}
\begin{proof}
The variational principle becomes
\begin{eqnarray*}
0 & = & \sum_{n=1}^N\Bigg(
\left\langle \dede{l}{\mathrm{X}}(\mathrm{X}^{n-1})+
(\mathrm{Q}^{n-1})^\dagger\mathrm{P}^{n-1},\delta\mathrm{X}^{n-1}
\right\rangle + 
\\
& & \qquad \left\langle
\delta\mathrm{P}^{n-1},\frac{1}{\Delta t}
\mathrm{Q}^{n-1}\logapp
\left(\left(\mathrm{Q}^{n-1}\right)^{-1}\mathrm{Q}^n\right)
-\mathrm{Q}^{n-1}\mathrm{X}^{n-1}
\right\rangle + \\
& & \qquad \Bigg\langle \mathrm{P}^{n-1},
\frac{1}{\Delta t}\mathrm{Q}^{n-1}
\left(T_{\left(\mathrm{Q}^{n-1}\right)^{-1}\mathrm{Q}^n}
\logapp\right)
\left(\mathrm{Q}^{n-1}\right)^{-1} \\
& & \qquad\qquad\qquad\qquad
\left(\delta\mathrm{Q}^n - \delta\mathrm{Q}^{n-1}\left(\mathrm{Q}^{n-1}\right)^{-1}\mathrm{Q}^n
\right)
+ \\
 & & \qquad\qquad \qquad\qquad
\delta\mathrm{Q}^{n-1}\left(\frac{1}{\Delta t}
\logapp\left(\left(\mathrm{Q}^{n-1}\right)^{-1}\mathrm{Q}^n
\right)
-\mathrm{X}^{n-1}\right)
\Bigg\rangle\Bigg).
\end{eqnarray*}
The discrete Clebsch equations are then
\begin{eqnarray}
0 & = & 
\dede{l}{\mathrm{X}}(\mathrm{X}^{n-1}) - (\mathrm{Q}^{n-1})^\dagger\mathrm{P}^{n-1}
, \\
0 & = & \mathrm{Q}^{n-1}\logapp
\left(\left(\mathrm{Q}^{n-1}\right)^{-1}\mathrm{Q}^n\right)
-\Delta t \mathrm{Q}^{n-1}\mathrm{X}^{n-1}
, \label{disc Q} \\
\nonumber 
0 & = & 
\left(\left(\mathrm{Q}^{n-1}\right)^{-1}\right)^\dagger
\left(T_{\left(\mathrm{Q}^{n-1}\right)^{-1}\mathrm{Q}^n}\logapp\right)^\dagger
\left(\mathrm{Q}^{n-1}\right)^\dagger\mathrm{P}^{n-1}
- \\ 
\nonumber & & \quad 
\left(\left(\mathrm{Q}^n\right)^{-1}\right)^\dagger
\left(T_{\left(\mathrm{Q}^n\right)^{-1}\mathrm{Q}^{n+1}}\logapp\right)^\dagger
\left(\mathrm{Q}^n\right)^\dagger \\
 & & \qquad \qquad \qquad
\mathrm{P}^n\left(\left(\mathrm{Q}^n\right)^{-1}\mathrm{Q}^{n+1}\right)
+ 
\\
& & \quad \left(
\frac{1}{\Delta t}\logapp\left(
\left(\mathrm{Q}^n\right)^{-1}\mathrm{Q}^{n+1}\right)\mathrm{Q}-
\mathrm{X}^n\right)^\dagger\mathrm{P}^n.
\end{eqnarray}
After rearrangement, equation \eqref{disc Q} becomes
\[
\mathrm{Q}^{n+1} = \mathrm{Q}^n\expapp({\Delta t}\mathrm{X}^n),
\]
and the last equation simplifies to
\begin{eqnarray*}
0 & = & \left(\left(\mathrm{Q}^{n-1}\right)^{-1}\mathrm{Q}^n\right)^\dagger
\left(T_{\left(\mathrm{Q}^{n-1}\right)^{-1}\mathrm{Q}^n}\logapp\right)^\dagger
\left(\mathrm{Q}^{n-1}\right)^\dagger\mathrm{P}^{n-1} + \\
& & \quad
\left(T_{\left(\mathrm{Q}^n\right)^{-1}\mathrm{Q}^{n+1}}\logapp\right)^\dagger
\left(\mathrm{Q}^n\right)^\dagger\mathrm{P}^n\left(\left(\mathrm{Q}^n\right)^{-1}\mathrm{Q}^{n+1}\right), \\
& = & \left(\expapp(\Delta t\mathrm{X}^{n-1})\right)^\dagger
\left(T_{\expapp(\Delta t\mathrm{X}^{n-1})}\logapp\right)^\dagger
\left(\mathrm{Q}^{n-1}\right)^\dagger\mathrm{P}^{n-1} - \\
& & \quad
\left(T_{\expapp(\Delta t\mathrm{X}^n)}\logapp\right)^\dagger
\left(\mathrm{Q}^n\right)^\dagger\mathrm{P}^n\left(\expapp(\Delta t\mathrm{X}^n)\right), 
\end{eqnarray*}
as required in the statement of Theorem \ref{discrete Clebsch}.
\qed\end{proof}
\begin{corollary}
$\mathrm{P}$ and $\mathrm{Q}$ can be eliminated from the equations arising from
the discrete variational principle in Definition \ref{euler} to 
obtain the equation
\begin{eqnarray*}
0 & = &  \left(\expapp(\Delta t\mathrm{X}^{n-1})\right)^\dagger
\left(T_{\expapp(\Delta t\mathrm{X}^{n-1})}\logapp\right)^\dagger
\dede{l}{\mathrm{X}}(\mathrm{X}^{n-1}) - \\
& & \qquad \qquad \left(T_{\expapp(\Delta t\mathrm{X}^n)}\logapp\right)^\dagger
\dede{l}{\mathrm{X}}(\mathrm{X}^n)
\left(\expapp(\Delta t\mathrm{X}^n)\right).
\end{eqnarray*}
\end{corollary}
\begin{proof}
  Substitute equation \eqnref{discrete dede l} into equation
  \eqnref{discrete Q}.  \qed\end{proof}

\begin{figure}
\begin{center}
\includegraphics[width=9cm]{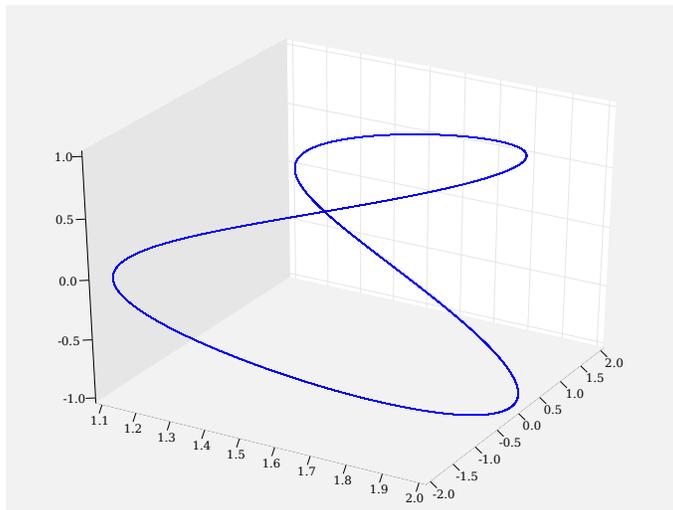}
\end{center}
\caption{\label{rigid plot} A plot showing numerical integration
over 100 periodic orbits using the discrete Clebsch integrator for
the rotating rigid body with moment of inertia eigenvalues (0.5,0.6,1)
and $\Delta t=0.1$. The good energy conservation and exact angular
momentum conservation are illustrated through the persistence of the 
periodic orbit structure over this long time integration interval.}
\end{figure}

\paragraph{Cayley transform methods} In the following example, we
derive an integrator for the rigid body equations by approximating the
exponential map using the Cayley transform, which preserves the
property of mapping from the Lie algebra into the group. This property
of the Cayley transform has long been used for ensuring that numerical
schemes preserve Lie group structure. In \cite{AuKrWa1993}, the Cayley
transform was used to reconstruct the attitude of a rotating rigid
body from numerical solutions of the body angular momentum equation
obtained from the midpoint rule. It was noted that the conservation of
the Casimir $\|\MM{m}\|^2$ (where $\MM{m}$ is the angular momentum)
was necessary to obtain conservation of spatial angular
momentum. \cite{LeSi1994} proposed to transform a Hamiltonian system
on a Lie group onto the Lie algebra using the Cayley transform (rather
than the exponential map), and integrating the resulting equation
numerically. This approach was developed in \cite{Is2001} which
suggested that, rather than integrating the Lie algebra equation using
a Runge-Kutta method (producing a Cayley Munthe-Kaas method
\cite{MK1998}), one could obtain an efficient scheme by using a
truncated Magnus expansion and, if a suitable numerical quadrature is
used to approximate the integrals in the series, one obtains a
time-reversible method of even order. In the example below, we embed
the Cayley transform into the discrete Clebsch variational principle;
higher-order schemes could be produced by using the methods of
Munthe-Kaas and Iserles.
\begin{example}[Rigid body integrator]
  For the case where $V$ is $\mathfrak{so}(3)$ and acts on
  $SO(3)$, we may approximate the exponential map to first-order using
  the Cayley transform
  \[
  \expapp(\mathrm{X})=\left(I + \frac{\mathrm{X}}{2}\right) 
  \left(I - \frac{\mathrm{X}}{2}\right)^{-1}.
  \]
  We obtain a corresponding approximation to the logarithm by taking 
  the inverse:
  \begin{eqnarray*}
    \logapp A &=& \expapp^{-1} A , \\
    &=& 2(A-I)(A+I)^{-1} = (A-I)\left(I + \frac{A-I}{2}\right)^{-1},\\
    &=& (A-I) + \mathcal{O}\left((A-I)^2\right),\\
    &=& \logapp(A) + \mathcal{O}\left((A-I)^2\right).
  \end{eqnarray*}
  Our approximation to the time derivative is then
\begin{eqnarray*}
\phi_{\Delta t}(\mathrm{Q}^{n+1},\mathrm{Q}^n)
&=& \mathrm{Q}^n\logapp((\mathrm{Q}^n)^{-1}\mathrm{Q}^{n+1}), \\
& = & \left(\frac{\mathrm{Q}^{n+1}-\mathrm{Q}^n}{\Delta t}\right)
\left(I + \frac{1}{2}\left(
(\mathrm{Q}^n)^{-1}\mathrm{Q}^{n+1} - I\right)\right)^{-1},
\end{eqnarray*}
which is the usual linear difference with a projection
applied to ensure that $\phi_{\Delta t}$ is in $\mathfrak{so}(3)$.

The $\mathrm{Q}$- and $\mathrm{X}$-equations are then
\begin{eqnarray}
\label{rigid xi} \dede{l}{\mathrm{X}}(\mathrm{X}^n) &=& (\mathrm{Q}^n)^T\mathrm{P}^n, \\
\label{rigid Q}
\mathrm{Q}^{n+1}\left(I-\Delta t\frac{\mathrm{X}^n}{2}\right) & = & \mathrm{Q}^n
\left(I+\Delta t\frac{\mathrm{X}^n}{2}\right).
\end{eqnarray}
The $\mathrm{Q}$-component of the variational principle, which
gives rise to the $\mathrm{P}$-equation, is
\[
\sum_{n=1}^N
\left\langle
\mathrm{P}^{n-1},\mathrm{Q}^{n-1}\delta\logapp\left(
(\mathrm{Q}^{n-1})^{-1}\mathrm{Q}^n\right)
\right\rangle=0.
\]
Making use of the formula
\begin{eqnarray*}
\delta\left(\logapp \left((\mathrm{Q}^{n-1})^{-1}\mathrm{Q}^n\right)\right)
\left(\frac{\left((\mathrm{Q}^{n-1})^{-1}\mathrm{Q}^n\right) + I}{2}\right) 
+ & & \\
\quad \logapp \left((\mathrm{Q}^{n-1})^{-1}\mathrm{Q}^n\right)
\frac{\delta \left((\mathrm{Q}^{n-1})^{-1}\mathrm{Q}^n\right)}{2} &=& 
\delta \left((\mathrm{Q}^{n-1})^{-1}\mathrm{Q}^n\right),
\end{eqnarray*}
we have
\begin{eqnarray*}
\delta\left(\logapp \left((\mathrm{Q}^{n-1})^{-1}\mathrm{Q}^n\right)\right)
&=&\left(I - 
\frac{1}{2}\logapp \left((\mathrm{Q}^{n-1})^{-1}\mathrm{Q}^n\right)\right) \\
& & \qquad
\delta \left((\mathrm{Q}^{n-1})^{-1}\mathrm{Q}^n\right)
\left(\frac{\left((\mathrm{Q}^{n-1})^{-1}\mathrm{Q}^n\right) + I}{2}\right)
^{-1}, \\
&=&\left(\frac{\left((\mathrm{Q}^{n-1})^{-1}\mathrm{Q}^n\right)+I}{2}\right)^{-1}
\\
& & \qquad
\delta \left((\mathrm{Q}^{n-1})^{-1}\mathrm{Q}^n\right)
\left(\frac{\left((\mathrm{Q}^{n-1})^{-1}\mathrm{Q}^n\right) + I}{2}\right)
^{-1}
, \\
&=&\left(\frac{\left((\mathrm{Q}^{n-1})^{-1}\mathrm{Q}^n\right)+I}{2}\right)^{-1}
(\mathrm{Q}^{n-1})^{-1} \\
& & \qquad \left(
\delta\mathrm{Q}^n - 
\delta\mathrm{Q}^{n-1}(\mathrm{Q}^{n-1})^{-1}\mathrm{Q}^n
\right)
\left(\frac{\left((\mathrm{Q}^{n-1})^{-1}\mathrm{Q}^n\right) + I}{2}\right)
^{-1}.
\end{eqnarray*}
Consequently the $\mathrm{P}$-equation in this formulation is
\begin{eqnarray*}
 & & \left((Q^{n-1})^{-1}\right)^T
\left(\left(\frac{\left((\mathrm{Q}^{n-1})^{-1}\mathrm{Q}^n\right) + I}{2}\right)
^{-1}\right)^T \\
 & & \qquad(\mathrm{Q}^{n-1})^T\mathrm{P}^{n-1}
\left(\left(\frac{\left((\mathrm{Q}^{n-1})^{-1}\mathrm{Q}^n\right) + I}{2}\right)
^{-1}\right)^T
\\
& = &\left((Q^{n})^{-1}\right)^T
\left(\left(\frac{\left((\mathrm{Q}^{n})^{-1}\mathrm{Q}^{n+1}\right) 
+ I}{2}\right)
^{-1}\right)^T \\
 & & \qquad(\mathrm{Q}^{n})^T\mathrm{P}^{n}
\left((\mathrm{Q}^n)^{-1}\mathrm{Q}^{n+1}\right)^T
\left(\left(\frac{\left((\mathrm{Q}^{n})^{-1}\mathrm{Q}^{n+1}\right)
 + I}{2}\right)
^{-1}\right)^T.
\end{eqnarray*}
Making use of equations (\ref{rigid xi}-\ref{rigid Q}) allows this
to be rearranged into the more compact form,
\begin{eqnarray*}
 & & \left(\expapp(\Delta t\mathrm{X}^{n-1})\right)^T
\left(\left(\frac{I + \expapp(\Delta t\mathrm{X}^{n-1})}{2}
\right)^{-1}\right)^T \\
& & \qquad\qquad \dede{l}{\mathrm{X}}(\mathrm{X}^{n-1})
\left(
\left(\frac{I + \expapp(\Delta t\mathrm{X}^{n-1})}{2}
\right)^{-1}
\right)^T
\\
& = & 
\left(\left
(\frac{I + \expapp(\Delta t\mathrm{X}^n)}{2}
\right)^{-1}\right)^T
\dede{l}{\mathrm{X}}(\mathrm{X}^n)
\expapp(\Delta t\mathrm{X}^n)^T
\left(\left(\frac{I + \expapp(\Delta t\mathrm{X}^n)}{2}
\right)^{-1}\right)^T
\end{eqnarray*}
Finally, making use of the Cayley transform approximation,
\[
\expapp(\Delta t\mathrm{X}) = \left(I + \Delta t\frac{\mathrm{X}}{2}\right)
\left(I - \Delta t\frac{\mathrm{X}}{2}\right)^{-1},
\]
and its consequence,
\[
I + \expapp(\Delta t\mathrm{X})
= 2\left(I - \Delta t\frac{\mathrm{X}}{2}\right)^{-1},
\]
we obtain the reduced equation for $\mathrm{X}$:
\begin{eqnarray*}
 & & 
\left(I - \Delta t\frac{\mathrm{X}^{n-1}}{2}\right)
\dede{l}{\mathrm{X}}(\mathrm{X}^{n-1})
\left(I + \Delta t\frac{\mathrm{X}^{n-1}}{2}\right)
\\
& = & 
\left(I + \frac{\Delta t\mathrm{X}^n}{2}\right)
\dede{l}{\mathrm{X}}(\mathrm{X}^n)
\left(I - \Delta t\frac{\mathrm{X}^n}{2}\right).
\end{eqnarray*}
This finally rearranges to become
\begin{eqnarray*}
\dede{l}{\mathrm{X}}(\mathrm{X}^n) &=& \dede{l}{\mathrm{X}}(\mathrm{X}^{n-1})
+ \frac{\Delta t}{2}\left(
\dede{l}{\mathrm{X}}(\mathrm{X}^{n-1})\mathrm{X}^{n-1}
+\dede{l}{\mathrm{X}}(\mathrm{X}^n)\mathrm{X}^n
\right) \\
& & \qquad
+ \frac{\Delta t}{2}\left(
\left(\mathrm{X}^{n-1}\right)^T
\dede{l}{\mathrm{X}}(\mathrm{X}^{n-1})
+\left(\mathrm{X}^n\right)^T\dede{l}{\mathrm{X}}(\mathrm{X}^n)
\right) , \\
& & \quad
+\frac{(\Delta t)^2}{4}\left(
\mathrm{X}^n
\dede{l}{\mathrm{X}}(\mathrm{X}^n)
\mathrm{X}^n
-
\mathrm{X}^{n-1}
\dede{l}{\mathrm{X}}(\mathrm{X}^{n-1})
\mathrm{X}^{n-1}
\right)
\end{eqnarray*}
which is the discrete rigid body equation obtained from
this choice of discrete Clebsch variational principle.
\end{example}
\begin{remark}
  The Clebsch integrator for the rigid body obtained using the Cayley
  transform approximation for $\exp$ is equivalent to the
  CAY-integrator obtained from the discrete Hamilton-Pontryagin
  principle in \cite{BoMa2007}. In that case, the equations are
  obtained by extremising a functional on a Lie algebra (in this case
  the kinetic energy as a function of the body angular velocity
  $\mathrm{X}$) subject to the constraint that
  $\dot{\mathrm{Q}}\mathrm{Q}^{-1} = \mathrm{X}$; in the case of the
  CAY-integrator the exponential map is again discretised using the
  Cayley transform. The Hamilton-Pontryagin principle provides an
  alternative viewpoint to the Clebsch principle; the extra feature in
  the Clebsch framework is the role of the $\diamond$-operator as a
  momentum map.
\end{remark}
A plot of the dynamics obtained from this discrete integrator
is given in Figure \ref{rigid plot}.

\section{Summary and outlook}
This paper has discussed Clebsch variational principles from the point
of view of a velocity map $\MM{\dot{Q}} = \,
\mathcal{L}_{\MM{\xi}}\MM{Q}$ which allows the dynamics $\MM{Q}(t)$ on
a manifold to be controlled by a time-series $\MM{\xi}(t)$ of elements
of a vector space, using Lagrange multipliers $\MM{P}(t)$. Theorem
\ref{reduction} shows that $\MM{Q}$ and $\MM{P}$ may be eliminated
from the resulting Clebsch equations, if and only if the velocity map
$\mathcal{L}_{\MM{\xi}}$ is a Lie algebra action on $\MM{Q}\in M$. The
Clebsch framework for velocity maps thus has a clear connection with
the theory of Euler-Poincar\'e reduction; namely, the equations
obtained are the Euler-Poincar\'e equations on the dual of the Lie
algebra. For the continuous time case where the velocity map is
assumed to be a Lie derivative, this connection is not unexpected.

Examples in the paper included the finite-dimensional rigid body
equation, and two infinite-dimensional examples: the singular
solutions of the EPDiff equation and the incompressible Euler
equation. In the EPDiff example, the Clebsch framework derives the
singular solutions as a family of momentum maps.

Finally the paper showed how discrete Clebsch variational principles
can be used to produce numerical methods for Clebsch equations.  For
the case of finite-dimensional Lie groups, in which the variational
principle need only be discretised in time, one may again eliminate
$\MM{Q}$ and $\MM{P}$ using the discrete approximation for the time
derivative in Definition \ref{Discretetimederiv} to obtain a
conservative numerical method in terms of $\MM{\xi}$ only. The example
of discretisation for the rigid body resulted in the CAY integrator
for the associated EP equation in \cite{BoMa2007}. Possible extensions
of this technique would be to obtain higher-order time-integration
methods based on Runge-Kutta/Munthe-Kaas methods \cite{MK1998} or
Magnus methods \cite{Is2001}, and to apply the Clebsch integrator to
systems with potentials such as the heavy top.

In the case of infinite-dimensional systems, it is necessary to
discretise the variational principle in space as well as time.  if one
wishes to eliminate $\MM{Q}$ and $\MM{P}$ thereby obtaining a closed
discrete equation for $\MM{\xi}$, Theorem \eqnref{reduction} requires
the spatial discretisation of the velocity map to remain a Lie algebra
action on the discretised space. This is a key step for making future
progress in applying the Clebsch method for discretising fluid
dynamics and other infinite-dimensional evolutionary systems. In the
special case where the Lagrangian is at most linear in space-time
derivatives (without higher derivatives), then the resulting PDE is
multisymplectic (see \cite{BrRe2001} and cited papers). Any discrete
Clebsch variational principle leads to a multisymplectic Clebsch
integrator. See \cite{CoHoHy2007} for more details.

\bibliography{Clebsch}
\nocite{*}
\end{document}